\documentclass[12pt]{article} 
\usepackage{latexsym,amssymb,amsmath}
  \newtheorem{theorem}{Theorem}
  \newtheorem{proposition}{Proposition}
               
               \newtheorem{lemma}{Lemma}
               
               \newtheorem{remark}{Remark}
               
               \def\pf{\par\noindent {\em Proof.}~\par\noindent}
               \def\qed{~\hfill{$\square$}\pagebreak[1]\par\medskip\par}

\newcommand{\f}{{\bf{f}}}
\newcommand{\Li}{{\mbox{Lip}}}

\newcommand{\cW}{{\cal W}}

\newcommand{\mD}{{\mbox D}}

\newcommand{\R}{{\mathbb R}}

\newcommand{\La}{\cL_{\alpha,\beta}}

\newcommand{\poz}{\partial_{\overline{z}}}
\newcommand{\pz}{\partial_{z}}

\newcommand{\cC}{{\cal C}}
\newcommand{\cQ}{{\cal Q}}
\newcommand{\cL}{{\cal L}}

\newcommand{\cT}{{\cal T}}

\begin{document}
\title{On the plane Lam\'e-Navier system in fractal domains}
\author
{Diego Esteban Gutierrez Valencia$^{(1)}$; Ricardo Abreu Blaya$^{(1)}$;\\Mart\'in Patricio \'Arciga Alejandre$^{(1)}$; Arsenio Moreno Garc\'ia$^{(2)}$}
\vskip 1truecm
\date{\small $^{(1)}$ Facultad de Matem\'aticas, Universidad Aut\'onoma de Guerrero, M\'exico.\\$^{(2)}$ Facultad de Inform\'atica y Matem\'atica, Universidad de Holgu\'in, Cuba\\Emails: diegogutierrez@uagro.mx, rabreublaya@yahoo.es,\\mparciga@gmail.com, amorenog@uho.edu.cu}
\maketitle
\begin{abstract}
\noindent
This paper is devoted to study a fundamental system of equations in plane Linear Elasticity Theory, the two-dimensional Lam\'e-Navier system. We rewrite them in a compressed form in terms of the Cauchy-Riemann operators and it allows us to solve a kind of Riemann problem for this  system. A generalized Teodorescu operator, to be introduced here, provides the means for obtaining the explicit solution of this problem for a very wide classes of regions, including those with a fractal boundary.  
\end{abstract}

\vspace{0.3cm}

\small{
\noindent
\textbf{Keywords.} Cauchy integral, fractals, linear elasticity.\\
\noindent
\textbf{Mathematics Subject Classification (2020).} 30G35.}
\section{Introduction}
The displacement equations of an elastic body under the  body force $(X,Y)$ in the case of plane strain reduce to
\begin{eqnarray}\label{lame1}
(\lambda + \mu)\frac{\partial \theta}{\partial x}+\mu \Delta u + X  =   0 \\
(\lambda + \mu)\frac{\partial \theta}{\partial y}+\mu \Delta v + Y  =   0\nonumber
\end{eqnarray}
where $\lambda$ and $\mu$ are the so-called Lam\'e parameters; $\Delta$ stands for the Laplacian in $\R^2$ and
\[
\theta=\frac{\partial u}{\partial{x}}+\frac{\partial v}{\partial y}.
\]
In the absence of body forces, i.e., when $X=Y=0$, the equations \eqref{lame1}  becomes
\begin{eqnarray}\label{lame0}
(\lambda + \mu)\frac{\partial \theta}{\partial x}+\mu \Delta u=0 \label{1.8}\\[0.1cm]
(\lambda + \mu)\frac{\partial \theta}{\partial y}+\mu \Delta v= 0 \label{1.9}\nonumber
\end{eqnarray}
As is standard, throughout the paper we will require that $\mu>0$ and $\lambda >-\frac{2\mu}{3}$.

Complex analysis has long been used to solve problems of plane elasticity. In view of the well-known Kolosov-Muskhelishvili formula (\cite{Mu})
the displacement field (as well as the stress function) can be represented by two holomorphic functions. Closely related to this formula is the fact (see \cite{Go}) that any bi-harmonic function $u$, i.e. $\Delta\Delta u=0$, in a simply connected domain $\Omega\subset\R^2$, is the real part of a bi-analytic complex function $f$. The Kolosov-Muskhelishvili formula is used in \cite{Mu} to study the Dirichlet problem of finding the elastic equilibrium of a body, if the displacements of the points on its boundary are known. The solutions of this problem follow from integral functional equations, which are obtained by using the Cauchy integral. However, to look for an explicit solution of this fundamental boundary value problem in elasticity theory for arbitrary plane domains is far from being trivial.

In this paper we rewrite the system \eqref{lame1} in terms of the Cauchy-Riemann operator $\poz=\frac{1}{2}(\frac{\partial}{\partial x}+i\frac{\partial}{\partial y})$ and its complex conjugate $\pz=\frac{1}{2}(\frac{\partial}{\partial x}-i\frac{\partial}{\partial y})$, which yields a sort of factorization of that system. Then, and without requiring the use of the Kolosov-Muskhelishvili  formula, we prove a Borel-Pompeiu integral representation  for $C^2$ functions. Properly formulated in our context this integral representation leads to a Cauchy representation formula for the solutions of \eqref{lame0} as well as to a particular solution of the inhomogeneous system \eqref{lame1}. 

Our approach is similar in spirit to the iterated method of Begehr \cite{Be} and goes back to Teodorescu \cite{Te}. Finally, and no less important, our results are extended to a more complicated framework of plane domains with fractal boundaries.
\section{Preliminaries}
In the sequel, unless stated otherwise, our attention will be concentrated on the system of the plane theory of elasticity when the region $\Omega\subset\R^2$ occupied by the body is assumed to be simply connected with a smooth closed Jordan curve $\gamma$ as its boundary. If necessary, we shall use the temporary notation $\Omega_+=\Omega$, $\Omega_-=\R^2\setminus (\Omega\cup\gamma)$.

As classically, consider the complex valued function $f=u+iv$ in the variable $z=x+iy$. After some simple transformations in \eqref{lame1} (with the entry of $i$ in this real-valued scenario), we obtain
\[
(\lambda+\mu)(\partial_{x}\theta+i\partial_{y}\theta)+\mu[\Delta u+i\Delta v]=-X-iY,
\]
or equivalently
\[
2(\lambda+\mu)\poz\theta+\mu\Delta f=-X-iY.
\]
But we know that $4\poz\pz=\Delta $. On the other hand, we have $\theta=\pz f+\overline{\pz f}$, which gives $\poz\theta=\poz\pz f+\poz\poz\bar{f}$. 

Consequently, we arrive to the following complex form of \eqref{lame1}:
\begin{equation*}
\left( \frac{\mu+\lambda}{2}\right) \partial_{\overline{z}}\partial_{\overline{z}}\overline{f}+\left(\frac{3\mu+\lambda}{2} \right) \partial_{\overline{z}}\partial_{z}f=g(z),
\end{equation*}
where $g=-\frac{1}{2}(X+iY)$.

In the absence of body forces, the above system becomes the homogeneous one
\[
\left( \frac{\mu+\lambda}{2}\right)\partial_{\overline{z}}\partial_{\overline{z}}\overline{f}+\left(\frac{3\mu+\lambda}{2} \right)\partial_{\overline{z}}\partial_{z}f=0,
\]
which in terms of the Poisson ratio $\sigma=\frac{\lambda}{2(\lambda+\mu)}$ reads
\begin{equation}\label{LMCp}
\partial_{\overline{z}}\partial_{\overline{z}}\overline{f}+(3-4\sigma)\partial_{\overline{z}}\partial_{z}f=0.
\end{equation}
An important point to note here is that the above complex reformulation straightforwardly suggests the form of the so-called {\emph{universal displacements}} of \eqref{lame0}: those displacement fields that can be maintained by applying boundary tractions in the absence of body forces for any linear elastic body, whatever the values of $\mu$ and $\lambda$.

In fact, from \eqref{LMCp} it follows that $f=u+iv$ is a universal solution if and only if in $\Omega$     
\[
\partial_{\overline{z}}\partial_{\overline{z}}\overline{f}=0,\,\partial_{\overline{z}}\partial_{z}f=0.
\]
Then $\pz f$ is simultaneously holomorphic and anti-holomorphic in $\Omega$ and hence a complex constant there. So we have $\pz f=A$, which yields
\[f(z)=Az+\overline{\varphi(z)},\,\,z\in\Omega,\]
where $A\in\mathbb{C}$ and $\varphi$ is an arbitrary holomorphic function in $\Omega$. 

In what follows, from consideration of symmetry and for the sake of brevity we put $\alpha=\frac{\mu+\lambda}{2}$ and $\beta=\frac{3\mu+\lambda}{2}$, giving rise to a notationally simpler complex Lam\'e-Navier system
\begin{equation}
\alpha\partial_{\overline{z}}\partial_{\overline{z}}\overline{f}+\beta\partial_{\overline{z}}\partial_{z}f=g(z).
\label{LMC1}
\end{equation}
The second order partial differential operator appearing in the left-hand side of \eqref{LMC1} will be denoted by $\cL_{\alpha,\beta}$, i.e. 
\[
\cL_{\alpha,\beta}[f]=\alpha\partial_{\overline{z}}\partial_{\overline{z}}\overline{f}+\beta\partial_{\overline{z}}\partial_{z}f.
\] 
\section{Borel-Pompeiu formula }
The fundamental theorem in our context is a Borel-Pompeiu formula, properly formulated to involve the Lam\'e-Navier operator $\La$. For further use, we introduce the notation 
\[
\alpha^*=\frac{\alpha}{\alpha^2-\beta^2},\,\,\beta^*=\frac{\beta}{\alpha^2-\beta^2}.
\]
\begin{theorem}[Borel-Pompeiu]\label{T1}
Let be $\Omega$ a simple connected domain with smooth boundary $\gamma$ and let $f\in C^2(\Omega)\cap C^1(\overline{\Omega})$. Then for $z\in\Omega$ we have
\begin{eqnarray}\label{bp}
f(z)=-\frac{\alpha\alpha^{*}}{2\pi i}{\int\limits_\gamma\frac{{f}(\xi)}{\overline{\xi-z}}\overline{d\xi}}-\frac{\beta \beta^{*}}{2\pi i}\int\limits_\gamma\frac{f(\xi)}{\xi-z}d\xi+\nonumber\\
\frac{\alpha^{*}}{2\pi i}\int\limits_{\gamma}\frac{{\xi-z}}{\overline{\xi-z}}\bigg[\alpha\partial_{\xi}f\overline{d\xi}+\beta\partial_{\overline{\xi}}\bar{f}\,\overline{d\xi}\bigg]+\nonumber\\
\frac{\beta^{*}}{2\pi i}\int\limits_{\gamma}\ln\left|\xi - z \right|^2\bigg[\alpha\partial_{\overline{\xi}}\bar{f}\,d\xi-\beta\partial_{\overline{\xi}}{f}\,\overline{d\xi}\bigg]+\nonumber\\
\frac{1}{\pi}\int\limits_{\Omega}\bigg[\alpha^{*}\frac{\xi-z}{\overline{\xi-z}}\overline{\mathcal{L}_{\alpha,\beta}[f(\xi)]}-{\beta^{*}}\ln\left|\xi-z \right|^2 \mathcal{L}_{\alpha,\beta}[f(\xi)]\bigg]d\xi.
\end{eqnarray}
\end{theorem}
\pf The proof consists in the convenient use of four already known identities. 

On one hand
\[
\frac{1}{\pi}\int_{\Omega}\frac{\xi-z}{\overline{\xi-z}}\partial_{\xi}\partial_{\overline{\xi}}\overline{f}d\xi=-\frac{1}{2\pi i}\int_{\gamma}\frac{\xi-z}{\overline{\xi-z}}\partial_{\overline{\xi}}\overline{f} \, \overline{d\xi}-\frac{1}{\pi}\int_{\Omega}\frac{1}{\overline{\xi-z}}\partial_{\overline{\xi}}\overline{f}d\xi
\]
and
\[
\frac{2}{\pi}\int_{\Omega}\ln\left|\xi-z \right| \partial_{\overline{\xi}}\partial_{\overline{\xi}}\overline{f}d\xi=\frac{1}{\pi i}\int_{\gamma}\ln\left|\xi-z \right| \partial_{\overline{\xi}}\overline{f} \, d\xi-\frac{1}{\pi}\int_{\Omega}\frac{1}{\overline{\xi-z}}\partial_{\overline{\xi}}\overline{f}d\xi,
\]
which follow directly from the Gauss theorem, on the other hand, the integral representation formulas
\begin{equation}\label{B1}
f(z)=-\frac{1}{2\pi i}\int_{\gamma}\frac{f(\xi)}{\overline{\xi-z}}\overline{d\xi}+\frac{1}{2\pi i}\int_{\gamma}\frac{\xi-z}{\overline{\xi-z}}\partial_{\xi}f(\xi)\overline{d\xi}+\frac{1}{\pi}\int_{\Omega}\frac{\xi-z}{\overline{\xi-z}}\partial_{\xi}\partial_{\xi}f(\xi)d\xi,
\end{equation}
\begin{equation}\label{B2}
f(z)=\frac{1}{2\pi i}\int_{\gamma}\frac{f(\xi)}{\xi-z}d\xi+\frac{1}{\pi i}\int_{\gamma}\ln \left|\xi-z \right| \partial_{\overline{\xi}}f(\xi)\overline{d\xi}+\frac{1}{\pi}\int_{\Omega}\ln\left|\xi-z \right|^2 \partial_{\xi}\partial_{\overline{\xi}}f(\xi)d\xi,
\end{equation}
obtained by H. Begehr \cite[p.227-228]{Be}.

Having disposed of this preliminary step, we rearranged the above formulas in order to obtain
\begin{equation*}
\frac{\alpha^{*}}{\pi}\int_{\Omega}\frac{\xi-z}{\overline{\xi-z}}\overline{\mathcal{L}_{\alpha,\beta}[f(\xi)]}d\xi-\frac{\beta^{*}}{\pi}\int_{\Omega}\ln\left|\xi-z \right|^2 \mathcal{L}_{\alpha,\beta}[f(\xi)]d\xi
\end{equation*}
\begin{equation*}
=\alpha\alpha^{*}\left\lbrace f(z)+\frac{1}{2\pi i} \int_{\gamma}\frac{1}{\overline{\xi-z}}f(\xi)\overline{d\xi}-\frac{1}{2\pi i}\int_{\gamma} \frac{\xi-z}{\overline{\xi-z}}\partial_{\xi}f(\xi)\overline{d\xi}\right\rbrace 
\end{equation*}
\begin{equation*}
+\beta\alpha^{*}\left\lbrace-\frac{1}{2\pi i}\int_{\gamma}\frac{\xi-z}{\overline{\xi-z}}\partial_{\overline{\xi}}\overline{f(\xi)} \, \overline{d\xi}-\frac{1}{\pi}\int_{\Omega}\frac{1}{\overline{\xi-z}}\partial_{\overline{\xi}}\overline{f(\xi)}d\xi \right\rbrace 
\end{equation*}
\begin{equation*}
-\alpha\beta^{*}\left\lbrace\frac{1}{2\pi i}\int_{\gamma}\ln\left|\xi-z \right|^2 \partial_{\overline{\xi}}\overline{f(\xi)} \, d\xi - \frac{1}{\pi}\int_{\Omega}\frac{1}{\overline{\xi-z}}\partial_{\overline{\xi}}\overline{f(\xi)}d\xi \right\rbrace 
\end{equation*}
\begin{equation*}
-\beta\beta^{*}\left\lbrace f(z)-\frac{1}{2\pi i}\int_{\gamma}\frac{1}{\xi-z}f(\xi) \, d\xi-\frac{1}{2\pi i}\int_{\gamma}\ln\left|\xi-z \right|^2 \partial_{\overline{\xi}}f(\xi) \, \overline{d\xi}  \right\rbrace.
\end{equation*}
Taking into account the obvious identities
\begin{equation}\label{identities}
\alpha\alpha^{*}-\beta\beta^{*}=1,\,\,\beta\alpha^{*}+\alpha\beta^{*}=0
\end{equation}
and after rather simple algebraic transformations we arrive at the desired relation \eqref{bp}.\qed
\begin{remark}
It should be noted that formulas \eqref{bp}, \eqref{B1} and \eqref{B2} remain valid in $\Omega_-$ if in the left-hand side of such formulas $f(z)$ is replaced by $0$. The proof of this standard fact is completely analogous (and even easier) to the case $z\in\Omega$.
\end{remark}

Inspired in the above formula we introduce the Teodorescu type operator 
\begin{equation}\label{Teo}
\cT_\Omega^{\cL}[g](z)=\frac{1}{\pi}\int\limits_{\Omega}\bigg[\alpha^{*}\frac{\xi-z}{\overline{\xi-z}}\overline{g(\xi)}-{\beta^{*}}\ln\left|\xi-z \right|^2 g(\xi)\bigg]d\xi,
\end{equation}
which runs as a right inverse for the Lam\'e-Navier operator $\La$. Indeed, we have
\begin{proposition}\label{inverse}
Let $g\in C(\overline{\Omega})$, then we have 
\[
\La\cT_\Omega^{\cL}[g](z)=\biggl\{
\begin{array}{rl}
g(z),\,\,& z\in\Omega_+\\
0,\,\,& z\in\Omega_-.
\end{array} 
\]
\end{proposition}
 \pf  We shall restrict ourselves to the case $z\in\Omega_+$. The case $z\in\Omega_-$ follows using similar
considerations.

After applying the operator $\mathcal{L}_{\alpha,\beta}$ to $\cT_\Omega^{\cL}[g(z)]$ we obtain:
\begin{equation}\label{tu1}
\mathcal{L}_{\alpha,\beta}\left\lbrace \cT_\Omega^{\cL}[g(z)]\right\rbrace=\alpha\partial_{\overline{z}}\partial_{\overline{z}}\overline{\cT_\Omega^{\cL}[g(z)]}+\beta\partial_{\overline{z}}\partial_{z}\cT_\Omega^{\cL}[g(z)].
\end{equation}
A direct calculation on the first term in the right-hand side of \eqref{tu1} gives
\begin{eqnarray}\label{tu2}
\alpha\partial_{\overline{z}}\partial_{\overline{z}}\overline{\cT_\Omega^{\cL}[g(z)]}&=&\frac{\alpha\alpha^{*}}{\pi}\partial_{\overline{z}}\partial_{\overline{z}}\left[ \int_{\Omega}\frac{\xi-z}{\overline{\xi-z}}g(\xi) \, d\xi\right]\!-\!\frac{\alpha\beta^{*}}{\pi}\partial_{\overline{z}}\partial_{\overline{z}}\left[ \int_{\Omega}\ln\left|\xi-z \right|^2 \overline{g(\xi)}d\xi\right] \nonumber	\\
&=&-\frac{\alpha\alpha^{*}}{\pi}\partial_{\overline{z}}\int_{\Omega}\frac{g(\xi)}{\xi-z} \, d\xi+\frac{\alpha\beta^{*}}{\pi}\partial_{\overline{z}}\int_{\Omega}\frac{\overline{g(\xi)}}{\overline{\xi-z}}d\xi\nonumber\\&=&\alpha\alpha^{*}g(z)+\frac{\alpha\beta^{*}}{\pi}\partial_{\overline{z}}\int_{\Omega}\frac{\overline{g(\xi)}}{\overline{\xi-z}}d\xi. 
\end{eqnarray}
Here we used the well-known fact that (see \cite[p. 31]{Ve}) 
\[
\partial_{\overline{z}}\bigg[-\frac{1}{\pi}\int_{\Omega}\frac{g(\xi)}{\xi-z}d\xi\bigg]=g(z)\,\, \mbox{in}\,\, \Omega.
\]
By a similar argument
\begin{equation}\label{tu3}
\beta\partial_{\overline{z}}\partial_{z}\cT_\Omega^{\cL}[g(z)]=-\beta\beta^{*}g(z)-\frac{\beta\alpha^{*}}{\pi}\partial_{\overline{z}}\int_{\Omega}\frac{\overline{g(\xi)}}{\overline{\xi-z}}d\xi
\end{equation}
Adding \eqref{tu2} and \eqref{tu3} we obtain
\begin{eqnarray*}
\mathcal{L}_{\alpha,\beta}\left\lbrace \cT_\Omega^{\cL}[g(z)]\right\rbrace & =  & \alpha\alpha^{*}g(z)+\frac{\alpha\beta^{*}}{\pi}\partial_{\overline{z}}\int_{\Omega}\frac{\overline{g(\xi)}}{\overline{\xi-z}}d\xi\\
& -  &\beta\beta^{*}g(z)-\frac{\beta\alpha^{*}}{\pi}\partial_{\overline{z}}\int_{\Omega}\frac{\overline{g(\xi)}}{\overline{\xi-z}}d\xi  \\
& =  & [\alpha\alpha^{*}-\beta\beta^{*}]g(z) = g(z),
\end{eqnarray*}
which is clear from \eqref{identities}.\qed

Of course, the above proposition means that the displacement given by $\cT_\Omega^{\cL}[g]$ represents a particular solution of the system
\eqref{LMC1} in presence of the the body force $g$.

An important corollary of the Theorem \ref{T1} is the following Cauchy representation formula for the elements of $\ker\La$.
\begin{theorem}
Let be $\Omega$ a simple connected domain with smooth boundary $\gamma$ and let $f\in C^2(\Omega)\cap C^1(\overline{\Omega})$. If moreover $f$ satisfies the homogeneous system \eqref{LMC1} (with $g\equiv 0$), then in $\Omega$ it can be represented by
\begin{eqnarray}\label{cf}
f(z)=-\frac{\alpha\alpha^{*}}{2\pi i}{\int\limits_\gamma\frac{{f}(\xi)}{\overline{\xi-z}}\overline{d\xi}}-\frac{\beta \beta^{*}}{2\pi i}\int\limits_\gamma\frac{f(\xi)}{\xi-z}d\xi+\nonumber\\
\frac{\alpha^{*}}{2\pi i}\int\limits_{\gamma}\frac{{\xi-z}}{\overline{\xi-z}}\bigg[\alpha\partial_{\xi}f\overline{d\xi}+\beta\partial_{\overline{\xi}}\bar{f}\,\overline{d\xi}\bigg]+\nonumber\\
\frac{\beta^{*}}{2\pi i}\int\limits_{\gamma}\ln\left|\xi - z \right|^2\bigg[\alpha\partial_{\overline{\xi}}\bar{f}\,d\xi-\beta\partial_{\overline{\xi}}{f}\,\overline{d\xi}\bigg].
\end{eqnarray}
\end{theorem}
The above Cauchy formula says that if the solution of the Dirichlet problem 
\begin{equation}\label{DP}
\Biggl\{
\begin{array}{rl}
\cL_{\alpha,\beta}F=0\,\,\mbox{in}\,\,\Omega\\
F=f\,\,\mbox{in}\,\,\gamma,
\end{array} 
\end{equation}
satisfying all the imposed conditions ($F\in C^2(\Omega)\cap C^1(\overline{\Omega})$) exists, it is necessarily given by \eqref{cf}. Of course, this solution is unique in accordance with the \emph{Uniqueness Theorem} \cite{Mu}, whose earlier proof as far as we know goes back to G. Kirchhoff \cite{Kir}.
\section{The Lam\'e-Cauchy transform}
Now let us introduce in our context a corresponding Cauchy transform specifically related to the system $\La f=0$, whose formal expression is suggested by \eqref{cf}. 

A suitable function space for this transform would be the space of $(1+\nu)$-order Lipschitz functions on $\gamma$.

We follow \cite{DAB1} (see also \cite{DAB2}) in defining $\Li(1+\nu,\gamma)$ ($0<\nu<1$) as the space of collections (Whitney jets) 
\[
\f:=\{f_0,f_1,f_2\}
\] 
of bounded functions defined in $\gamma$ and such that for all $t,\tau\in\gamma$
\begin{eqnarray*}
&&|f_0(t)-f_0(\tau)-(t-\tau){f}_1(\tau)-(\overline{t-\tau}){f}_2(\tau)|\leq c|t-\tau|^{1+\nu}\\
&&|{f}_1(t)-{f}_1(\tau)|\leq c|t-\tau|^{\nu},\,\,|{f}_2(t)-{f}_2(\tau)|\leq c|t-\tau|^{\nu},
\end{eqnarray*}
the constant $c>0$ being independent of $t$, $\tau$.

At first glance, the introduction of the space $\Li(1+\nu,\gamma)$ seems to be artificial. However, this space is closely connected with a very deep result in Real Analysis, the so-called \emph{Whitney Extension Theorem}. Here we state without proof a complex version of this fundamental theorem
\begin{theorem}[Whitney]\label{Wh}
Let $\f$ be in $\Li(1+\nu,\gamma)$.  Then, there exists a compact supported complex valued function $\tilde{f}\in C^{1,\nu}(\R^2)$ satisfying
\begin{itemize}
\item[(i)] $\tilde{f}|_\gamma = f_0,\,\pz\tilde{f}|_\gamma=f_1,\poz\tilde{f}|_\gamma=f_2$
\item[(ii)] $\tilde{f}\in C^\infty(\R^{2} \setminus \gamma)$,
\item[(iii)] $|\partial_z^{j_1}\partial_{\overline{z}}^{j_2}\tilde{f}(z) | \leqslant c \, \mbox{\em dist}(z,\gamma)^{\nu-1}$, for $j_1+j_2=2$ and $z\in\R^2\setminus\gamma$.
\end{itemize}
\end{theorem}
The pioneering work here is \cite{Wh}. An excellent reference alone more classical lines is the well-known book of E. M. Stein \cite[Chapter VI, p. 176]{St}.

The already announced  Cauchy operator (the sum of contour integrals in \eqref{cf}) may be naturally defined for collections $\f\in\Li(1+\nu,\gamma)$. These collections of functions are intrinsically given on $\gamma$ and play a similar roll as H\"older functions do for the classical holomorphic Cauchy integral.

More precisely, for $\f\in\Li(1+\nu,\gamma)$ we define the Lam\'e-Cauchy transform by
\begin{eqnarray}\label{ct}
\cC^\cL{\f}(z)=-\frac{\alpha\alpha^{*}}{2\pi i}{\int\limits_\gamma\frac{{f}_0(\xi)}{\overline{\xi-z}}\overline{d\xi}}-\frac{\beta \beta^{*}}{2\pi i}\int\limits_\gamma\frac{f_0(\xi)}{\xi-z}d\xi+\nonumber\\
\frac{\alpha^{*}}{2\pi i}\int\limits_{\gamma}\frac{{\xi-z}}{\overline{\xi-z}}\bigg[\alpha f_1\overline{d\xi}+\beta\bar{f_1}\overline{d\xi}\bigg]+\nonumber\\
\frac{\beta^{*}}{2\pi i}\int\limits_{\gamma}\ln\left|\xi - z \right|^2\bigg[\alpha\bar{f_1}\,d\xi-\beta{f_2}\,\overline{d\xi}\bigg].
\end{eqnarray}

The behavior of $\cC^\cL{\f}(z)$ at infinity is obviously governed by the last integral. Indeed, the first two terms are in fact Cauchy type integrals vanishing at $z=\infty$, while the third term is bounded. Consequently, we have 
\begin{equation}\label{infty}
\cC^\cL{\f}(z)={\mathcal{O}}(\ln|z|),\,\,z\to\infty.
\end{equation}
 
As was to be expected, the function given by $\cC^\cL\f$ satisfies in $\R^2\setminus\gamma$ the homogeneous Lam\'e system \eqref{LMC1}.

\begin{theorem}
Let $\f\in\Li(1+\nu,\gamma)$, then
\[
\La[\cC^\cL\f]=0,\,\,\mbox{in}\,\,\Omega_+\cup\Omega_-.
\]
\end{theorem}  
 \pf For simplicity we omit the tedious but straightforward calculations involved in the first part of the proof. After that we get
\begin{eqnarray*}
\mathcal{L}_{\alpha,\beta}[\cC^\cL\f](z) & =  & \frac{\alpha \beta \beta^{*}}{\pi i}\int_{\gamma}\frac{\overline{f_0(\xi)}}{(\overline{\xi - z})^3}\overline{d\xi}\!-\!\frac{\beta^2\alpha^{*}}{2\pi i}\int_{\gamma}\frac{\overline{f_1(\xi)}}{(\overline{\xi-z})^2}\overline{d\xi}\!-\!\frac{\beta\alpha\alpha^{*}}{2\pi i}\int_{\gamma}\frac{f_1(\xi)}{(\overline{\xi-z})^2}\overline{d\xi} \\
	& +  &\frac{\alpha^2\beta^{*}}{2\pi i}\int_{\gamma}\frac{f_1(\xi)}{(\overline{\xi-z})^2} \, \overline{d\xi}-\frac{\alpha\beta\beta^{*}}{2\pi i}\int_{\gamma}\frac{\overline{f_2(\xi)}}{(\overline{\xi-z})^2} \, d\xi,
\end{eqnarray*}
and hence 
\begin{equation}\label{after}
\mathcal{L}_{\alpha,\beta}[\cC^\cL\f](z)=\frac{\alpha \beta \beta^{*}}{\pi i}\!\!\int_{\gamma}\frac{\overline{f_0(\xi)}}{(\overline{\xi - z})^3} \overline{d\xi}\!-\!\frac{\beta^2\alpha^{*}}{2\pi i}\!\int_{\gamma}\frac{\overline{f_1(\xi)}}{(\overline{\xi-z})^2}\overline{d\xi}\!-\!\frac{\alpha\beta\beta^{*}}{2\pi i}\!\int_{\gamma}\frac{\overline{f_2(\xi)}}{(\overline{\xi-z})^2}d\xi,
\end{equation}
since $\alpha \beta^{*}=\beta\alpha^{*}$.

In the next step use will be made of the Whitney extension theorem. Indeed, let be $\tilde{f}\in C^{1,\nu}(\R^2)$ as in Theorem \ref{Wh},  then \eqref{after} becomes
\[
\mathcal{L}_{\alpha,\beta}[\cC^\cL\f](z)=\frac{\alpha \beta \beta^{*}}{\pi i}\int_{\gamma}\frac{\overline{\tilde{f}(\xi)}}{(\overline{\xi - z})^3}\overline{d\xi}-\frac{\beta^2\alpha^{*}}{2\pi i}\int_{\gamma}\frac{\partial_{\overline{\xi}}\overline{\tilde{f}(\xi)}}{(\overline{\xi-z})^2}\overline{d\xi}-\frac{\alpha\beta\beta^{*}}{2\pi i}\int_{\gamma}\frac{\partial_{\xi}\overline{\tilde{f}(\xi)}}{(\overline{\xi-z})^2}d\xi.
\]
Assume $z\in\Omega_-$, then in accordance with
\[
\partial_{\overline{\xi}}\left(\frac{\overline{\tilde{f}(\xi)}}{\overline{(\xi-z)}^2} \right)=\frac{\partial_{\overline{\xi}}\overline{\tilde{f}(\xi)}}{\overline{(\xi-z)}^2}-2\frac{\overline{\tilde{f}(\xi)}}{\overline{(\xi-z)}^3},
\]
the above formula reads
\begin{eqnarray}\label{after1}
\mathcal{L}_{\alpha,\beta}[\cC^\cL\f](z) & =  & \frac{\beta^2\alpha^{*}}{2\pi i}\int_{\gamma}\frac{\partial_{\overline{\xi}}\overline{\tilde{f}(\xi)}}{(\overline{\xi-z})^2}\overline{d\xi}-\frac{\beta^2\alpha^{*}}{2\pi i}\int_{\gamma}\partial_{\overline{\xi}}\left(\frac{\overline{\tilde{f}(\xi)}}{(\overline{\xi-z})^2} \right) \, \overline{d\xi}\nonumber\\
& -  & \frac{\beta^2\alpha^{*}}{2\pi i}\int_{\gamma}\frac{\partial_{\overline{\xi}}\overline{\tilde{f}(\xi)}}{(\overline{\xi-z})^2} \, \overline{d\xi}-\frac{\beta^2\alpha^{*}}{2\pi i}\int_{\gamma}\frac{\partial_{\xi}\overline{\tilde{f}(\xi)}}{(\overline{\xi-z})^2} \, d\xi \nonumber\\
& =  & -\frac{\beta^2\alpha^{*}}{2\pi i}\left[\int_{\gamma}\partial_{\overline{\xi}}\left(\frac{\overline{\tilde{f}(\xi)}}{(\overline{\xi-z})^2} \right)\overline{d\xi} +\int_{\gamma}\frac{\partial_{\xi}\overline{\tilde{f}(\xi)}}{(\overline{\xi-z})^2}d\xi\right] \nonumber\\
& =  &  -\frac{\beta^2\alpha^{*}}{2\pi i}\left[\int_{\gamma}\partial_{\overline{\xi}}\left(\frac{\overline{\tilde{f}(\xi)}}{(\overline{\xi-z})^2} \right) \, \overline{d\xi} +\int_{\gamma}\partial_{\xi}\left( \frac{\overline{\tilde{f}(\xi)}}{(\overline{\xi-z})^2}\right)d\xi\right].
\end{eqnarray}
Gauss formula applied to the last expression in the right-hand side of \eqref{after1} leads to $\mathcal{L}_{\alpha,\beta}[\cC^\cL\f](z)=0$ for $z\in\Omega_-$.

Now let $z\in\Omega$. Removing from $\Omega$ the ball $B(z,\epsilon)$ gives the domain $\Omega_\epsilon:=\Omega\setminus B(z,\epsilon)$ with
boundary $\partial\Omega_\epsilon=\gamma\cup(-C_\epsilon)$, where $C_\epsilon$ is the circumference with radius $\epsilon$ and center $z$. Now, we applied Gauss theorem once again in the right-hand side of \eqref{after1}, but this time to the domain $\Omega_\epsilon$. 

Then we have
\begin{eqnarray*}  
\int_{\gamma}\partial_{\overline{\xi}}\left(\frac{\overline{\tilde{f}(\xi)}}{(\overline{\xi-z})^2} \right) \, \overline{d\xi} +\int_{\gamma}\partial_{\xi}\left( \frac{\overline{\tilde{f}(\xi)}}{(\overline{\xi-z})^2}\right)d\xi\\=\int_{C_\epsilon}\partial_{\overline{\xi}}\left(\frac{\overline{\tilde{f}(\xi)}}{(\overline{\xi-z})^2} \right) \, \overline{d\xi} +\int_{C_\epsilon}\partial_{\xi}\left( \frac{\overline{\tilde{f}(\xi)}}{(\overline{\xi-z})^2}\right)d\xi.
\end{eqnarray*}
But for $\xi\in C_\epsilon$ we have 
\[
\frac{1}{(\overline{\xi-z})^2}=\frac{(\xi-z)^2}{\epsilon^4}
\]
and hence
\begin{eqnarray*}  
\int_{\gamma}\partial_{\overline{\xi}}\left(\frac{\overline{\tilde{f}(\xi)}}{(\overline{\xi-z})^2} \right) \, \overline{d\xi} +\int_{\gamma}\partial_{\xi}\left( \frac{\overline{\tilde{f}(\xi)}}{(\overline{\xi-z})^2}\right)d\xi\\=\frac{1}{\epsilon^4}\int_{C_\epsilon}\partial_{\overline{\xi}}\left({{({\xi-z})^2}\overline{\tilde{f}(\xi)}} \right) \, \overline{d\xi} +\int_{C_\epsilon}\partial_{\xi}\left({({\xi-z})^2}{\overline{\tilde{f}(\xi)}}\right)d\xi=0,
\end{eqnarray*}
the last equality being a consequence of the Gauss theorem applied to $B(z,\epsilon)$.\qed

It should be noted that the weakly singularity of the kernels $\frac{{\xi-z}}{\overline{\xi-z}}$ and $\ln\left|\xi - z \right|^2$ implies that the third and fourth integrals in \eqref{ct} do not experience a jump when $z$ is crossing the boundary $\gamma$. This together with the classical Plemelj-Sokhotski formulas applied to the first two Cauchy type integrals in \eqref{ct} lead to
\[
[\cC^\cL{\f}]^+(t)-[\cC^\cL{\f}]^-(t)=f_0(t),
\]
where \[[\cC^\cL{\f}]^\pm (t)=\lim_{\Omega_\pm\ni z\to t}[\cC^\cL{\f}](z).\]
On the other hand, it is also possible to deal with the jump of the function $\pz[\cC^\cL{\f}(z)]$ on the boundary $\gamma$. For this purpose, we use Theorem \ref{Wh} and rewrite \eqref{ct} as
\begin{eqnarray}\label{ctW}
\cC^\cL{\f}(z)={\alpha\alpha^{*}}\bigg[-\frac{1}{2\pi i}{\int\limits_\gamma\frac{\tilde{f}(\xi)}{\overline{\xi-z}}\overline{d\xi}}+\frac{1}{2\pi i}\int\limits_{\gamma}\frac{{\xi-z}}{\overline{\xi-z}}\partial_\xi\tilde{f}(\xi)\overline{d\xi}\bigg]\nonumber\\-\beta \beta^{*}\bigg[\frac{1}{2\pi i}\int\limits_\gamma\frac{\tilde{f}(\xi)}{\xi-z}d\xi+\frac{1}{2\pi i}\int\limits_{\gamma}\ln\left|\xi - z \right|^2\partial_{\overline{\xi}}\tilde{f}(\xi)\,\overline{d\xi}\bigg]\nonumber\\+
\alpha^*\beta\bigg[  
\frac{1}{2\pi i}\int\limits_{\gamma}\frac{{\xi-z}}{\overline{\xi-z}}\partial_{\overline{\xi}}\overline{\tilde{f}}(\xi)\overline{d\xi}+\frac{1}{2\pi i}\int\limits_{\gamma}\ln\left|\xi - z \right|^2 \partial_{\overline{\xi}}\overline{\tilde{f}}(\xi)d\xi\bigg]
\end{eqnarray}
After using the identities \eqref{B1} and \eqref{B2} in the first two expressions in brackets, we obtain
\begin{eqnarray}\label{ctW+}
\cC^\cL{\f}(z)={\alpha\alpha^{*}}\bigg[\tilde{f}(z)-\frac{1}{\pi}\int_{\Omega}\frac{\xi-z}{\overline{\xi-z}}\partial_{\xi}\partial_{\xi}\tilde{f}(\xi)d\xi\bigg]\nonumber\\-\beta \beta^{*}\bigg[ \tilde{f}(z)-\frac{1}{\pi}\int_{\Omega}\ln\left|\xi-z \right|^2 \partial_{\xi}\partial_{\overline{\xi}}f(\xi)d\xi\bigg]\nonumber\\+
\alpha^*\beta\bigg[  
\frac{1}{2\pi i}\int\limits_{\gamma}\frac{{\xi-z}}{\overline{\xi-z}}\partial_{\overline{\xi}}\overline{\tilde{f}}(\xi)\overline{d\xi}+\frac{1}{2\pi i}\int\limits_{\gamma}\ln\left|\xi - z \right|^2 \partial_{\overline{\xi}}\overline{\tilde{f}}(\xi)d\xi\bigg]
\end{eqnarray}
for $z\in\Omega_+$.

Then in $\Omega_+$ we have

\begin{eqnarray}\label{D(ctW+)}
\pz\cC^\cL{\f}(z)={\alpha\alpha^{*}}\bigg[\pz\tilde{f}(z)+\frac{1}{\pi}\int_{\Omega}\frac{1}{\overline{\xi-z}}\partial_{\xi}\partial_{\xi}\tilde{f}(\xi)d\xi\bigg]\nonumber\\-\beta \beta^{*}\bigg[ \pz\tilde{f}(z)+\frac{1}{\pi}\int_{\Omega}\frac{1}{\xi-z}\partial_{\xi}\partial_{\overline{\xi}}\tilde{f}(\xi)d\xi\bigg]\nonumber\\+
\alpha^*\beta\bigg[  
-\frac{1}{2\pi i}\int\limits_{\gamma}\frac{1}{\overline{\xi-z}}\partial_{\overline{\xi}}\overline{\tilde{f}}(\xi)\overline{d\xi}-\frac{1}{2\pi i}\int\limits_{\gamma}\frac{1}{\xi-z} \partial_{\overline{\xi}}\overline{\tilde{f}}(\xi)d\xi\bigg]
\end{eqnarray}
or equivalently
\begin{eqnarray}\label{D(ctW+1)}
\pz\cC^\cL{\f}(z)={\alpha\alpha^{*}}\bigg[-\frac{1}{2\pi i}\int\limits_{\gamma}\frac{1}{\overline{\xi-z}}\partial_{{\xi}}\tilde{f}(\xi)\overline{d\xi}\bigg]\nonumber\\-\beta \beta^{*}\bigg[\frac{1}{2\pi i}\int\limits_{\gamma}\frac{1}{{\xi-z}}\partial_{{\xi}}\tilde{f}(\xi){d\xi}\bigg]\nonumber\\+
\alpha^*\beta\bigg[  
\overline{\frac{1}{2\pi i}\int\limits_{\gamma}\frac{1}{{\xi-z}}\partial_{{\xi}}{\tilde{f}}(\xi){d\xi}}-\frac{1}{2\pi i}\int\limits_{\gamma}\frac{1}{\xi-z} \partial_{\overline{\xi}}\overline{\tilde{f}}(\xi)d\xi\bigg]
\end{eqnarray}

In an analogous manner one has for $z\in\Omega_-$
\begin{eqnarray}\label{D(ctW-1)}
\pz\cC^\cL{\f}(z)={\alpha\alpha^{*}}\bigg[-\frac{1}{2\pi i}\int\limits_{\gamma}\frac{1}{\overline{\xi-z}}\partial_{{\xi}}\tilde{f}(\xi)\overline{d\xi}\bigg]\nonumber\\-\beta \beta^{*}\bigg[\frac{1}{2\pi i}\int\limits_{\gamma}\frac{1}{{\xi-z}}\partial_{{\xi}}\tilde{f}(\xi){d\xi}\bigg]\nonumber\\+
\alpha^*\beta\bigg[\overline{\frac{1}{2\pi i}\int\limits_{\gamma}\frac{1}{{\xi-z}}\partial_{{\xi}}{\tilde{f}}(\xi){d\xi}}-\frac{1}{2\pi i}\int\limits_{\gamma}\frac{1}{\xi-z} \partial_{\overline{\xi}}\overline{\tilde{f}}(\xi)d\xi\bigg]
\end{eqnarray}

A careful look to the two last Cauchy type integrals in brackets reveals that they experience the same jump crossing $\gamma$ and then the jump of the difference itself vanishes. 

Consequently, we obtain
\[
[\pz\cC^\cL{\f}]^+(t)-[\pz\cC^\cL{\f}]^-(t)=(\alpha\alpha^*-\beta\beta^*)\pz\tilde{f}(t)=f_1.
\]  
We have thus proved that, given a Whitney jet $\f=\{f_0,f_1,f_2\}$ on $\gamma$, the function $F(z)=\cC^\cL{\f}(z)$ is a solution of the boundary value problem
\begin{equation}\label{JP}
\begin{array}{rl}
\La F(z)&=0,\,\,\,\,z\in\Omega_+\cup\Omega_-\\
F^+(t)-F^-(t)&=f_0(t),\,\,\,\,t\in\gamma,\\
{[\pz F]}^+(t)-{[\pz F]}^-(t)&=f_1(t),\,\,\,\,t\in\gamma
\end{array} 
\end{equation}
If we additionally assume that $\pz F(\infty)=0$ and $F(z)={\mathcal{O}}(\ln|z|)$ as $z\to\infty$, then  $F(z)=\cC^\cL{\f}(z)$ is the unique solution of \eqref{JP} up to an arbitrary complex constant. More precisely:
\begin{theorem}
Let $\f\in\Li(1+\nu,\gamma)$. The boundary value problem \eqref{JP} has a solution given by $F(z)=\cC^\cL{\f}(z)$. Moreover, it is unique (up to an arbitrary constant) under the asymptotic conditions $\pz F(\infty)=0$ and $F(z)={\mathcal{O}}(\ln|z|)$ as $z\to\infty$.
\end{theorem}
\pf The existence part was already proved. It remains to prove the uniqueness (up to a complex constant), which is carried out as usual indirectly. The assumption of two solutions $F_1, F_2$ implies that $F=F_1-F_2$ fulfills
\begin{equation}\label{JP0}
\begin{array}{rl}
\La F(z)&=0,\,\,\,\,z\in\Omega_+\cup\Omega_-\\
F^+(t)-F^-(t)&=0,\,\,\,\,t\in\gamma,\\
{[\pz F]}^+(t)-{[\pz F]}^-(t)&=0,\,\,\,\,t\in\gamma.
\end{array} 
\end{equation}
Since $\La F(z)=0$ in $\Omega_+\cup\Omega_-$, the function $G(z)=\alpha\overline{\pz F}+\beta\pz F$ is holomorphic there and has the jump 
\[
G^+(t)-G^-(t)=\alpha\overline{\{[{\pz F}]^+(t)-[{\pz F}]^-(t)\}}+\beta\{[{\pz F}]^+(t)-[{\pz F}]^-(t)\}=0,
\]
by the second boundary condition in \eqref{JP0}.

This together with the vanishing condition $G(\infty)=0$ yields $G\equiv 0$ and therefore we get $\pz F\equiv 0$, since $\alpha\not=\pm\beta$. 

Consequently $\overline{F(z)}$ is holomorphic in $\Omega_+\cup\Omega_-$ and has no jump through $\gamma$, which is clear from the first boundary condition in \eqref{JP0}. In this way we finally obtain that $\overline{F(z)}$ is holomorphic in the whole plane and it grows at a sufficiently slow rate ($\overline{F(z)}={\mathcal{O}}(\ln|z|)$), which reduces $\overline{F(z)}$ and hence $F(z)$ to a constant, by a slight generalization of the Liouville theorem.\qed  
\section{The fractal case} 
The solutions of boundary value problems in plane elasticity theory for regions of general shape presents considerable practical and theoretical difficulties. The main new ingredient of this section is the extension of our previous considerations to the case of domains $\Omega$ admitting a fractal boundary. Our method goes back to the early work of Boris Kats \cite{K}  in connection with Riemann boundary value problems for holomorphic functions (see also \cite{ABK}).

We follow \cite{HN} in assuming that $\gamma$ is $d$-summable for some $1<d<2$, i.e. the improper integral
\[
\int_0^1 N_\gamma(\tau) \, \tau^{d-1} \, d\tau
\]
converges, where $ N_\gamma(\tau)$ stands for the minimal number of balls of radius $\tau$ needed to cover $\gamma$.

As was early remarked in \cite{HN} any curve $\gamma$ with fractal box dimension $\mD(\gamma)$ is $d$-summable for any $d=\mD(\gamma)+\epsilon$, $\epsilon>0$. In particular, if $\gamma$ is so smooth as in the previous sections, then $\mD(\gamma)=1$ and hence it is $(1+\epsilon)$-summable for any $\epsilon>0$.

For the rest of the paper we assume $\Omega$ to be a Jordan domain with $d$-summable boundary $\gamma$.

The following lemma appeared in \cite[Lemma 2]{HN} and reveals the specific importance of the notion of $d$-summability applied to the boundary $\gamma$ of a Jordan domain $\Omega$ in connection with the Whitney decomposition $\cW$ of $\Omega$ by $k$-squares $\cQ$ of diameter $|\cQ|$.
\begin{lemma}\cite{HN}
\label{dsum}
If $\Omega$ is a Jordan domain of $\R^{2}$ and its boundary $\gamma$ is $d$-summable, then the expression $\sum_{\cQ\in{\mathcal{W}}}|\cQ|^d$, called the $d$-sum of the Whitney decomposition ${\cW}$ of $\Omega$, is finite.
\end{lemma}

\begin{lemma}\label{lp}
Let $\f\in\Li(1+\nu,\gamma)$, then $\partial_z^{j_1}\partial_{\overline{z}}^{j_2}\tilde{f}(z)\in L^p(\Omega)$, $j_1+j_2=2$  for $p=\frac{2-d}{1-\nu}$.
\end{lemma}
\pf In fact, we have
\begin{eqnarray*}
\int\limits_\Omega|\partial_z^{j_1}\partial_{\overline{z}}^{j_2}\tilde{f}(\xi)|^p d\xi=\hspace{-0.1cm}\sum_{\cQ\in{\cal W}}\int\limits_{\cQ}|\partial_z^{j_1}\partial_{\overline{z}}^{j_2}\tilde{f}(\xi)|^p d\xi\le c\,\hspace{-0.1cm}\sum_{\cQ\in{\cal W}}\int\limits_{\cQ}({\mbox{dist}}(\xi,\gamma))^{p(\nu-1)}d\xi,
\end{eqnarray*}
which follows from Theorem \ref{Wh} (iii).

Consequently,
\begin{eqnarray*}
\int\limits_\Omega|\partial_z^{j_1}\partial_{\overline{z}}^{j_2}\tilde{f}(\xi)|^p d\xi\le c\sum_{\cQ\in{\cal W}} |\cQ|^{p(\nu-1)}|\cQ|^{2}=c\hspace{-0.1cm}\sum_{\cQ\in{\cal W}} |\cQ|^{d}<+\infty,
\end{eqnarray*}
by Lemma \ref{dsum}.\qed  

In addition, notice that for $\nu>\frac{d}{2}$ we have $p=\frac{2-d}{1-\nu}>2$. Therefore, under this condition we have $\partial_z^{j_1}\partial_{\overline{z}}^{j_2}\tilde{f}\in L^p(\Omega)$, for some $p>2$ and the same is obviously true for $\La(\tilde{f})$. 

Based in the following theorem the boundary value problem \eqref{JP} has a solution, applicable even to our fractal context. 
\begin{theorem}\label{casifinal}
Let $\f\in\Li(1+\nu,\gamma)$ and let $\gamma$ be $d$-summable with $\nu>\frac{d}{2}$. The boundary value problem \eqref{JP} has a solution given by 
\begin{equation}\label{frS}
F(z)=\chi_\Omega(z)\tilde{f}(z)-\cT_\Omega^{\cL}[\La(\tilde{f})(z)],\,\,z\in\Omega_+\cup\Omega_-,
\end{equation}
where $\chi_\Omega$ stands for the characteristic function of $\Omega$.
\end{theorem} 
\pf 
The fact that $F(z)$ satisfies $\La f=0$ in $\Omega_+\cup\Omega_-$ is a direct consequence of Proposition \ref{inverse}.

On the other hand, the first jump condition in \eqref{JP} easily follows from the weakly singularity of the kernels $\frac{{\xi-z}}{\overline{\xi-z}}$ and $\ln\left|\xi - z \right|^2$, i.e.,  $\cT_\Omega^{\cL}[\La(\tilde{f})(z)]$ do not experiences a jump when $z$ is crossing the boundary $\gamma$. 

To prove the second jump condition we use Theorem 6.1 in \cite{Ve}. It is immediately clear that
\[
\pz F(z)=\chi_\Omega(z)\pz\tilde{f}(z)+\alpha^*\frac{1}{\pi}\overline{\int\limits_\Omega\frac{1}{{\xi-z}}{\La(\tilde{f})(\xi)}}d\xi-\beta^*\frac{1}{\pi}\int\limits_\Omega\frac{1}{\xi-z}\La(\tilde{f})(\xi)d\xi.
\] 
Since $\nu>\frac{d}{2}$, $\La(\tilde{f})$ is in $L^p(\Omega)$, with $p>2$, and then
\[
\frac{1}{\pi}\int\limits_\Omega\frac{1}{\xi-z}\La(\tilde{f})(\xi)d\xi
\]  
represents a continuous function in $\R^2$ and so the second jump condition in \eqref{JP} follows.\qed  
\begin{remark}
Working in domains with  fractal boundaries, the question of uniqueness of \eqref{frS} becomes more subtle since a fractal curve is not in general a removable set for continuous holomorphic functions. Following similar ideas to those used by Kats in \cite{K} this situation could be partially overcome, but we will not develop this point here.
\end{remark}

\end{document}